\documentclass[11pt]{article}

\usepackage{amsmath}
\usepackage{amsfonts}
\usepackage{amssymb}

\date{DATE: 17/10/2013}

\newtheorem{Lemma} {{\textsc Lemma}}
\newtheorem{Prop} {\textsc{Proposition}}

\newcommand{\QED}{\hfill $\square$\bigskip}
\newcommand{\qed}{\hfill $\square$\bigskip}
\newcommand{\pf}{{\noindent\bf Proof.\ \ }}
\newcommand{\N}{{\Bbb N}}
\newcommand{\Z}{{\Bbb Z}}
\newcommand{\Q}{{\Bbb Q}}
\newcommand{\J}{{\cal J}}
\newcommand{\la}{{\langle}}
\newcommand{\ra}{{\rangle}}
\newcommand{\e}{{\epsilon}}%
\newcommand{\G}{{\Gamma}}%

\newcommand{\iso}{\simeq}
\newcommand{\p}{{\varphi}}
\newcommand{\g}{{\gamma}}%

\begin{document}

\title{Inertial endomorphisms of an abelian group}

\author{Ulderico Dardano\footnote{Ulderico Dardano,  Dipartimento di Matematica e
Applicazioni ``R.Caccioppoli'', Universit\`a di Napoli ``Federico
II'', Via Cintia - Monte S. Angelo, I-80126 Napoli, Italy.
\texttt{dardano@unina.it}}\  \ -\
Silvana Rinauro\footnote{Silvana Rinauro,
Dipartimento di Matematica, Informatica ed Economia, Universit\`a della
Basilicata, Via dell'Ateneo Lucano 10 - Contrada Macchia Romana,
I-85100 Potenza, Italy. \texttt{silvana.rinauro@unibas.it} } }

\maketitle
\vskip-2cm 
\begin{abstract}
 We describe inertial endomorphisms of an abelian group $A$, that is endomorphisms $\varphi$  with the property $|(\varphi(X)+X)/X|<\infty$ for each $X\le A$. They form a ring containing multiplications, the so-called finitary endomorphisms and non-trivial instances.
 We show that  inertial invertible endomorphisms form a group, provided $A$ has
finite torsion-free rank.  In any case, the group $IAut(A)$ they generate
 is commutative modulo the group $FAut(A)$ of finitary automorphisms,
which is known to be locally finite. We deduce that $IAut(A)$ is locally-(center-by-finite). Also
we consider the lattice dual property, that is $|X/(X\cap \varphi(X))|<\infty$ for each $X\le A$. We show that this implies the above one, provided $A$ has finite torsion-free rank.
\end{abstract} 
\vskip-0.1cm 
\noindent{\bf 2000 Mathematics Subject Classification}. \\ 
\emph{Primary 20K30, Secondary 20E07, 20E36, 20F24}. 


\section{Introduction and statement of main results}

Recently there has been interest for totally inert (TIN) groups, i.e. groups whose all subgroups are inert (see \cite{B}, \cite{DGMT}, \cite{DET}, \cite{R}).
A subgroup is said \emph{inert} if it is commensurable to each conjugate of its.
 Two subgroup $X,Y$ of any group are told \emph{commensurable} iff
 $X\cap Y$ has finite index in both $X$ and $Y$ (see \cite{SH}).

When dealing with soluble TIN-subgroups of  groups one is concerned with automorphisms with the following property. As in \cite{DGBSV2} and \cite{DGBSV}, an \emph{endomorphism} $\p$ of an abelian group $A$ (from now on always in additive notation) is
said (right-) \emph{inertial} iff:\\
(RIN)\hskip2cm $\forall X\le A \ \ |\varphi(X)+X\ /X|<\infty$.\\
Consideration of endomorphisms instead of automorphisms is due to the
following.

\medskip \noindent
{\bf Fact} \emph{Inertial endomorphisms of any abelian group $A$ form a ring, say
$IE(A)$, containing the ideal $F(A)$ of endomorphisms  with finite image.}

\medskip
\noindent To prove this notice that  if $\p_1$ and $\p_2$ both have RIN, then $\forall X\le A$ \ \ $|(X+\p_1(X))/X|<\infty$ and $|(X +\p_1(X)+\p_2(X)+\p_1\p_2(X))/(X+\p_1(X))|<\infty$.

\noindent  In Theorem A below we give a rather satisfactory characterization of inertial endomorphisms of an abelian group, from which we deduce useful consequences. In Proposition A we exhibit non-trivial instances of inertial endomorphisms.

\medskip
\noindent\textbf{Corollary A} \emph{The ring $IE(A)/F(A)$ is commutative.}

\medskip
In \cite{DR} we considered inertial \emph{automorphisms} of an abelian group $A$  generalizing previous results from \cite{Cas} and \cite{FGN}. From results in \cite{DR} it follows that
for an automorphism $\p$ of a periodic abelian group $A$ above (RIN) is equivalent to its lattice dual condition (\emph{left-inertial}):\\
(LIN)\hskip2cm $\forall X\le A \ \ |X/(X\cap \varphi(X))|<\infty$.\\
Note that we will say just \emph{inertial} for RIN.
Endomorphisms which have both LIN and RIN are those  mapping subgroups to commensurable ones. Notice that $x\mapsto 2x$ is an automorphism of $\Q^\omega$ with RIN but not LIN.
Clearly, the inverse of an automorphism with RIN has LIN and vice-versa.

We
consider the group $IAut(A)$ generated by inertial automorphisms of $A$.
From Theorem A we have:

\medskip
\noindent\textbf{Corollary B}  \emph{Let $\p$ be an endomorphism of an abelian group $A$.\\
1)\ If $A$ has finite torsion-free rank, then LIN implies RIN and the two properties are equivalent if $\p$ is an  automorphism. Thus inertial automorphisms form the group $IAut(A)$.
\\
2) If $A$ has not finite torsion-free rank, then $IAut(A)$  is formed by the products $\g_1\g_2^{-1}$, where $\g_1,\g_2$ are both inertial automorphisms.
}

\medskip

Recall that the rank $r_0(A)$  of any free abelian subgroup $F$ of $A$ such that $A/F$ is periodic
is said \emph{torsion-free rank} of $A$. From Proposition \ref{PropMixENDO} we will have that \emph{when $r_0(A)=\infty$ an endomorphism is both LIN and RIN iff it acts as the identity or the inversion map on a subgroup with finite index}.

\medskip

\emph{Automorphisms acting as the identity map on a finite index subgroup form a group, say $FAut(A)$, which is locally finite} (see \cite{W1}) and \emph{is  contained in} $IAut(A)$. Actually, $IAut(A)$ is not periodic but we have a corresponding statement, which also follows from Theorem A.

\bigskip
\noindent\textbf{Theorem B} \emph{Let $\G=IAut(A)$ be the group generated by the inertial automorphisms
of an abelian group $A$. Then:\\
1)\ $\G'\le FAut(A)$ is locally finite;\\
2)\ $\G$ is locally central-by-finite.}\\

\medskip
\noindent
Thus periodic elements of $IAut(A)$ form a subgroup containing the derived subgroup. However, there are non-elementary instances of periodic non-finitary inertial automorphisms. To see this consider the $p$-group $B\oplus D$ where ($p\ne2$) $B$ is infinite bounded and $D$ is divisible with finite rank and the automorphism acting as the identity on $B$ and the inversion map on $D$. On the other side, the abelian group $IAut(A)/FAut(A)$ may be rather large as in next statement.

\bigskip
\noindent\textbf{Proposition A }  \emph{There exists a countable abelian group $A$ with $r_0(A)=1$ such that $IAut(A)$ has a subgroup $\Sigma\iso \prod_p \Z(p)$ with $ \Sigma\cap FAut(A)=T(\Sigma)\iso \bigoplus_p \Z(p)$, where $p$ ranges over the set of all primes.}

\medskip

We state now the main result of the paper. By multiplication of an abelian group $A$ we mean the (componentwise) action by $p$-adics, if the group is periodic; otherwise we mean the natural action by rational numbers (when possible and uniquely defined). \emph{Multiplications form a ring and commute with any  endomorhism}, clearly, and are often inertial, as in Proposition \ref{mult+iner+}.

\medskip
\noindent\textbf{Theorem A }
{\it Let $\p_1,\dots,\p_t$ finitely many endomorphisms of an abelian  group $A$. Then each
$\p_i$ is inertial if and only if there is a finite index subgroup $A_{0}$ of $A$ such that (a) or (b) holds:

\medskip\noindent
$(a)$ each $\p_i$ acts as multiplication by $m_i\in\Z$ on $A_{0}$;

\medskip\noindent
 $(b)$ $A_{0}=B\oplus D\oplus C$ and there exist finite sets of primes $\pi\subseteq \pi_1$ such that:

i)\ $B\oplus D$ is the $\pi_1$-component of $A_{0}$ where $B$ is bounded and $D$ is a divisible $\pi'$-group with finite rank,

ii)\ $C$ is a $\Q^{\pi}[\p_1,\dots,\p_t]$-module, with a submodule $V\iso \Q^{\pi}\oplus\dots\oplus\Q^{\pi}$ (finitely many times) such that $C/V$ is a $\pi_1$-divisible $\pi'$-group,

iii)\ each $\p_i$ acts by (possibly different) multiplications on $B$, $D$, $V$, $C/V$,

iv)\ each
$\p_i=\frac{m_i}{n_i}\in\Q$ on $V$ and on all $p$-components of $D$ such that the $p$-component of $C/V$ is infinite and $\pi=\pi(n_1\cdots n_t)$.
}

\bigskip\noindent
Notice that when $A$ is periodic, the above statement applies with  $V=0$ and amounts to Proposition \ref{TEOREMA_PERIODICO+} below.
On the other side, if $A$  is torsion-free, then inertial endomorphisms are just multiplications, see Proposition \ref{CasoTorsionFree++}. Further, for details on LIN condition see Propositions \ref{TEOREMA_PERIODICO+} and \ref{PropMixENDO}.


\section{Preliminaries and Terminology}

As a standard reference on abelian groups we use \cite{F}.
Letter $A$ always denote an abelian group, we regard as a \emph{left} $E(A)$-module, where $E(A)$ denotes the ring of endomorphism of $A$.
Letters $\p,\g$ denote endomorphisms of $A$, while $m,n,r,s,t$ denote integers, $p$ always a prime, $\pi$ a set of primes,
$\pi(n)$ the set of prime divisors of $n$. We denote: by $T=T(A)$ the torsion subgroup of $A$, by  $A_\pi$ the $\pi$-component of $A$, by $A[n]:=\{a\in A\ |\ na=0\ \}$. If $nA=0$, we say that $A$ is \emph{bounded} by $n$. Further, we say that $A$ is bounded, if it is bounded by some $n$. Denote by $D=Div(A)$ the largest divisible subgroup of $A$.

For a subset $X$ and an endomorphism $\p$ of a (left) $R$-module $A$, we denote
by $\la X\ra$  (resp. $X^{(\p)}=R[\p]X$) the additive subgroup (resp. the $R[\p]$-submodule) spanned by $X$, as usual. Note that, if $X$ is a $R$-submodule, then by $X _{(\p)}$ we mean the largest $R[\p]$-submodule contained in $X$.

When $A$ is $\pi$-divisible (that is $pA=A$ for each $p\in\pi$) and $A_\pi=0$, we regard it as a $\Q^{\pi}$-module and write $\p=\frac{m}{n}\in \Q^{\pi}$ for the well-defined map via $\p(nx)=mx$. As usual $\Q^{\pi}$ denote the ring
of rationals whose denominator is a $\pi$-number, where $\Q^{(n)}:=\Z[\frac{1}{n}]$.

We call \emph{multiplications} of an abelian group $A$ either the above actions of elements of $\frac{m}{n}\in\Q$ or, when $A$ is periodic, the componentwise actions
of the elements of the cartesian product of the rings $\Q_p^*$ of  $p$-adics for each prime $p$.
Recall that $\Q_p^*$ contains the ring of rationals whose denominator is coprime to $p$. Note that we use word ``multiplication" in a way different from \cite{F}. Ours are in fact ``scalar" multiplications.

When $A$ is periodic, multiplications are precisely the so-called \emph{power} endomorphisms  which are defined by $\forall X\le A \ \p(X)\subseteq X$ (see \cite{R}).
On the other hand, when $A$ is not periodic, power endomorphisms are multiplication $x\mapsto nx$ for fixed $n\in\Z$.

 Let us point out which multiplications have RIN or LIN. Recall that an
abelian group $A$ with the \emph{minimal condition} (Min) is just a group of the shape $A=F\oplus D$, where $F$ is finite and $D$ is divisible with finite total rank that is the sum of finitely many infinite cocyclic (Pr\"ufer) groups. For short we will write $A$ has \emph{FTFR} when the torsion-free rank
$r_0(A)$ is finite.

 \begin{Prop}\label{mult+iner+} Let $\p$ be a multiplication of an infinite abelian group $A$, then \\
\noindent R)\ \ $\p$ is inertial iff either $A$ has FTFR or $\p$ is a multiplication by an integer;\\
\noindent L1) if $A$ is a $p$-group, then $\p$ is LIN iff if $\p$ is invertible or $A$ has Min and $\p\ne 0$;\\
L2) if $0<r_0(A)<\infty$, then $\frac{m}{n}$ is LIN iff $A_{\pi(m)}$ has Min and $\p\ne 0$ ($m,n$ coprime) ;\\
L3) if $r_0(A)=\infty$, then $\p$ is LIN iff $\p=\frac{1}{n}$.\\
\end{Prop}

 \pf  Let us prove part (R). If $A$ is periodic, the statement is trivial.  Otherwise, let $\p=\frac{m}{n}$. Now if $\p$ is inertial, then, arguing in $\bar A:=A/T(A)$, we have that for any $\bar X$ free with infinite rank, it results that $(\p(\bar X)+\bar X)/\bar X$ is infinite.

Conversely, if $\p=m$ acts as an integer, it is trivially inertial. Assume $A$ has FTFR. For any $X\le A$, the section $(\p(X)+X)/X$ is a bounded $\pi(m)$-group and is finite mod $T:=T(A)$.
On the other hand, since
$A$ is a $\Q^{\pi(m)}$-module,  $A_{\pi(m)}=0$ and therefore $(\p(X)+X)/X$ avoids $T$.

 \medskip
 If $\p$ is LIN, then  $A/\p(A)<\infty$ implies that $\p\ne 0$. Further we have:

\medskip\noindent L1) Let the $p$-adic $\alpha=p^s\alpha_1$ represent $\p$ on $A$ with $\alpha_1$ invertible. If $\p$ is not invertible, then $s>0$ and $A[p]\le ker \p$ is finite. Hence $A$ has Min. Conversely, if $A$ has Min, then for any $X\le A$ we have that $X/\p(X)=X/p^sX$ is finite.

\medskip\noindent
L2)  If $\p=\frac{m}{n}$ then $A_{\pi(m)}$ has Min by $(L1)$. Conversely, for each $X\le A$ we have
$$\frac{X}{X\cap\frac{m}{n}X}\iso\frac{nX}{nX\cap mX}$$
if finite as bounded by $m$ and the rank of $A_{\pi(m)}$ and torsion-free rank of $A$ are finite.

\medskip\noindent
L3) Let $\p=\frac{m}{n}$ and take $\bar X\le \bar A:=A/T$ free of infinite rank. As above
 ${\bar X}/({\bar X\cap\frac{m}{n}\bar X)}$ is infinite unless $|m|=1$. On the other hand $\p=\frac{1}{n}$ is LIN as $X\le\p(X)$ for each $X\le A$.\qed

The other way round, let us see that inertial endomorphisms of a torsion-free abelian group
are all multiplications and there is a natural ring monomorphism $IE(A)\hookrightarrow \Q$.
\medskip 

\begin{Prop}\label{CasoTorsionFree++}
Let $\p$ be an endomorphism of a {torsion-free} abelian group$A$.\\
R) $\p$ is RIN iff $\p$ is multiplication by $\frac{m}{n}$ where if ${n}\ne\pm1$ then $r_0(A)<\infty$.\\
L) $\p$ is LIN iff $\p$ is multiplication by $\frac{m}{n}$ with $m\not=0$ and if $m\ne\pm1$ then $r_0(A)<\infty$.

\end{Prop}

\noindent In particular, if $\p\ne0$ and $r_0(A)<\infty$, $\p$ is LIN iff  $\p$ is RIN.
\noindent In the above statement $m$ and $n$ are meant to be coprime integers.

\pf The sufficiency of the condition is clear in both cases (R) and (L). To prove necessity, we generalize an argument used in \cite{DR}. In both cases (R) and (L), if $a\in A$  there exist  $m,n\in\Z$ such that
$ma=n\p(a)$. As $A$ is torsion-free, $m,n$ can be choosen coprime.
Let us show that $\frac{m}{n}$ is independent of
$a$. Let $a_1\in A$. If $\la a_{1} \ra\cap\la a
\ra\ne\{0\}$, then $ka_1=ha$ for some $h,k\in \Z$.
Therefore we may write $\p(a_1)=
\frac{h}{k}\p(a)=\frac{h}{k}\frac{m}{n}a=\frac{m}{n}a_1.$
Similarly, if $\la a_{1} \ra\cap\la a \ra=\{0\}$, there
will exist $m_1,m_2,n_1,n_2\in\Z$ such that
$$\frac{m}{n}a+\frac{m_1}{n_1}a_1=\p(a)+\p(a_1)=
\p(a+a_1)=
\frac{m_2}{n_2}(a+a_1)=\frac{m_2}{n_2}a+\frac{m_2}{n_2}a_1 .$$ It
follows $\frac{m}{n}=\frac{m_2}{n_2}=\frac{m_1}{n_1}$. Thus $\p$ is a multiplication.

For the rank restriction, apply Proposition \ref{mult+iner+}.
\qed

\bigskip

Note a sufficient condition for an endomorphism to be inertial. We omit the
straightforward proof.

\begin{Prop}\label{inerXfin} Let $\p$ be an endomorphism of an abelian  group $A$. If  $\p$ acts as an inertial (resp. LIN) endomorphism either on a finite index subgroup of $A$ or modulo a finite subgroup, then $\p$ is inertial (resp. LIN) on the whole $A$ indeed.\qed
\end{Prop}

We will use often the following fact.

\begin{Prop}\label{PF} For an endomorphism $\p$ of a periodic abelian group $A$ the following are equivalent:\\
MF)\ $\p$ acts as a multiplication on a finite index subgroup $A_0$ of $A$,\\
FM)\ $\p$ acts as a multiplication modulo a finite subgroup $A_1$ of $A$.
\end{Prop}

\noindent \pf This is very easy. If $(MF)$ holds and $\p=\alpha\in\J$ on $A_0$, where
$\J$ is the cartesian product of the rings $\Q_p^*$ of the $p$-adics for each prime $p$.
Then consider the natural action of $\alpha$ on $A$. Thus $A_1:=im(\p-\alpha)$ is an epimorphic image of $A/A_0$. The converse is similar. \qed


\section{Inertial endomorphisms of a \emph{periodic} abelian group}

This section is devoted to prove the periodic case of the main Theorem A of the paper. In fact we prove a more detailed characterization
of inertial endomorphisms of torsion abelian groups. We handle LIN-endomorphisms as well.

\begin{Prop}\label{TEOREMA_PERIODICO+}
Let  $\p_1,...\p_t$ be finitely many endomorphisms of an abelian periodic group $A$. Then:

\medskip
\noindent R)\ each $\p_i$ is inertial
 iff there is a finite index subgroup $A_0=B\oplus D\oplus C$ of $A$ such that:

i) $B\oplus D$ and $C$ are coprime,

ii) $B$ is bounded and $D$ is divisible with Min,

iii) each $\p_i$ is multiplication on $B$, $D$ and $C$.

\medskip
If the above holds and $\Phi:=\Z[\p_1,...\p_t]$, then

$(FS)\hskip 1cm \exists m\ \forall X\le A\ \ |X^{\Phi}/X_{\Phi}|\le m.$

 \medskip\noindent L)\ each $\p_i$ is LIN  iff it is inertial and there are subgroups $A_0,\ B,\ D,\ C$ as above such that $\p_i$ is a non-zero multiplication on each non-zero primary component of $D$ and an invertible multiplication on $B$ and $C$.
\end{Prop}

\noindent Recall that, for each $i$, $X^{\Phi}$ and $X_{\Phi}$ are $\p_i$-invariant and
$X^{\Phi}\le X\le X_{\Phi}$.

Before proving the Proposition we state the following easy but
fundamental fact, which reduces the proof to the case $A$ is a $p$-group.

\begin{Prop}\label{riduzione_ai_p-gruppi} An endomorphism
of an abelian torsion group $A$  is inertial (resp. LIN) iff it is such on all primary
components and multiplication (resp. invertible multiplication) on all but finitely many of them.\end{Prop}

\pf It is trivial that the condition is sufficient. Concerning necessity, we only deal with case  LIN, the case RIN being similar. Let $\pi$ be the set of primes $p$ such that $\p$ is not invertible multiplication on  $A_{p}$. If  $p\in \pi$, then  either $\p$ is not a multiplication on $A_{p}$ or $\p$ is a non-invertible multiplication. In the former case there is a cyclic subgroup $X_{p}$ of $A_{p}$ such that $\p(X_{p})\not\subseteq  X_{p}$, and hence $| X_{p}\cap \p(X_{p})|<|\p(X_{p})|\le|X_{p}|$. In the latter case there is a cyclic subgroup $X_{p}$ of $A_{p}$ such that $\p(X_{p})$ is properly contained in $X_{p}$. In both cases $|X_{p}/(X_{p}\cap \p(X_{p}))|>1$. It is now clear that if $\p$ is LIN, then $\pi$ is finite, as $|X/(X\cap \p(X))|$ must be finite for $X:=\bigoplus\limits_{p\in \pi}X_{p}$.
 \qed

\noindent
We prove a couple of Lemmas. The first one extends Proposition 4.3 in \cite{R}.

\begin{Lemma}\label{Robinson++} Let $A$ be an abelian $p$-group, $a\in A$ and $\p\in E(A)$.\\
1)\ \ If $\p$ either RIN or LIN, then the cyclic $\p$-submodule $\Z[\p](a)=:\la a\ra^{(\p)}$ of $A$ generated by $a$ is finite.\\
2)\ \ If $|X/X_{(\p)}|<\infty$ for all $X\le A$,
then $|X^{(\p)}/X|<\infty$ for all  $X\le A$.\\
3)\ \ If $|X/X_{(\p)}|\le p^{m}$ for all $X\le A$,
then $|X^{(\p)}/X|\le p^{m^{2}}$ for all  $X\le A$.
\end{Lemma}

\pf
(1) We may assume $A=\la a\ra^{(\p)}$. Suppose first $a$ has order prime $p$ and
 consider
the natural epimorphism  of $\Z_p[x]$-modules mapping $1$ to $a$ and $x$ to $\p(a)$ (regard $A$ as $\Z_p[x]$-module where $x$ acts as $\p$):
$$F:\Z_p[x]\mapsto A.$$
If $F$ is injective,
we can replace $A$ by $\Z_p[x]$ and $\p$ by multiplication by $x$. If
$H:=\Z_p[x^2]$, then $\p(H)=xH$ is infinite,
while $H\cap xH=0$, a contradiction. Then $F$ is not injective and $A$ is finite as it is isomorphic to a proper quotient of $\Z_p[x]$.
If now $a$ has (any) order $p^{\e}$, then $A/pA$ is finite, by the above. Moreover, $pA=\la pa\ra^{(\p)}$ is finite by induction on $\e$.

\medskip
(2) This can be proved in a similar way as case $(3)$

\medskip
(3) We claim that {\em if $a\in A$ has order $p^{\e}$, then
$|\la a\ra^{(\p)}|\le p^{(m+1)\e}.$}\\
Assume first  $\e=1$, that is $a$ has order $p$ and $A_0:=\la a\ra^{(\p)}$
is elementary abelian. Suppose, by contradiction, the above $F$
is injective. As above, let $H:=\Z_p[x^2]$. Then $H_{(\p)}=(g(x^2))$ for some polynomial $g$.
Since $|H/H_{(\p)}|=p^m<\infty$, we have $g\ne 0$. Then $(g(x^2))\not\subseteq H$, a contradiction.
Therefore, for some $f\in \Z_p[x]$ with degree say $n$, we have
$$\frac{\Z_p[x]}{(f)}\iso_\p \la a\ra^{(\p)}=A_0.$$
Thus the minimal $\p$-invariant subgroups of $A_0$ correspond $1-1$
to the irreducible monic factors of $f$, which are at most $n$. Consider a $\Z_p$-basis $X$ of A containing an element in each subgroup of them. The the hyperplane $H$ of equation
$x_{1}+x_{2}+\cdots+x_{n}=0$ has  index $p$ in  $\la a\ra^{(\p)}$ and  $H_{(\p)}=0$ as $H\cap X=\emptyset$. Therefore $ |\la a\ra^{(\p)}|\le p^{m+1}.$

If $\e>1$, by induction $B:=\la p^{\e-1}a \ra^{(\p)}$ has
order at most $p^{(m+1)(\e-1)}$ and $\la a\ra ^{(\g)}/B$ has
order at most $p^{m+1}$ by case $\e=1$. Therefore
$|a^{(\p)}|\le p^{(m+1)\e}$, as claimed.

In the general case let $X$ be any subgroup of $A$ and
$X_{(\p)}=0$. Thus $|X|=:p^{\e}\le p^{m}$.
Write
$X=\la a_1\ra \oplus\cdots\oplus \la a_r\ra$ with $a_{i}$ of order
$p^{\e_{i}}$ and  $\e_{1}+\cdots+\e_{r}= \e$. Since $|\la a_{i}\ra^{(\p)}|\le p^{(m+1)\e_{i}}$ by the above, we have
$|X^{(\p)}|\le p^{(m+1)\e}$. So that  $|X^{(\p)}/X|\le p^{(m+1)\e-\e}\le p^{m^{2}}$.
\qed

\begin{Lemma}\label{LemmaDiv}Let $D$ be a divisible periodic subgroup of an abelian subgroup $A$ and $\p\in E(A)$.
If $\p$ is either  RIN or LIN, then $\p$ is multiplication on $D$.
\end{Lemma}

\pf  Without loss of generality, we may assume $D$ and $\p(D)$ are Pr\"ufer groups. If $\p$ is LIN, then $D\le \p(D)$ and thus $D=\p(D)$. Therefore in both cases  RIN or LIN, we have $ \p(D)\le D$.\qed

\medskip
\noindent \textbf{Proof of Proposition \ref{TEOREMA_PERIODICO+}}
 Note that $(FS)$ implies trivially that each $\p_i$ is inertial. By Proposition \ref{riduzione_ai_p-gruppi} we may \emph{assume $A$ is a $p$-group} (see details below).
We  proceed by a sequence of claims and start by considering the case of a single endomorphism $\p:=\p_1$.

\medskip\noindent$\bullet$ \emph {\ If $A$ is elementary abelian and $\p$ is either inertial or LIN, then $\p$ is FM}. Assume by contradiction $\p$ is not FM (as in Proposition \ref{PF}). We prove that:
\emph{for any finite subgroup $X\le A$
 such that $X\cap \p(X)=0$ there exists a finite subgroup $X'> X$ such that} $$X'\cap \p(X')=0\ \ {\rm and } \ \ \p(X')>\p(X).$$  Therefore, starting at $X_0=0$, by recursion we
 define $X_{i+1}:=X_i'$ and $X_\omega:=\cup_i X_i$. We get that both $X_\omega$ and $\p(X_\omega)$ are
 infinite and $X_\omega\cap \p(X_\omega)=0$, a contradiction.

 First note that $Z:=X^{(\p)}$ is finite by Lemma \ref{Robinson++}.
 Since $\p$ does not act as a  multiplication on $A/Z$,
 we can choose $a\in
A$ such that
 $\p(a)\not\in \la a,Z\ra$ and define $X':=\la a\ra +X$.
If now
 $y\in X'\cap \p(X')$,  then $\exists m,n\in \Z_p$, $\exists x,x_0\in X$ such
 that $y=ma+x=n\p(a)+\p(x_0)$. Thus $n\p(a)\in \la a\ra+
Z$ while $\p(a)\not\in \la a\ra+Z$. Therefore $n=0$ and $ma=0$ as well. It follows $y=x=\p(x_0)\in X\cap
 \p(X)=0$, as claimed.


\medskip\noindent$\bullet$ \emph{If $\p$ is LIN then it  is inertial, as in both cases it holds}:
$$(fs) \hskip 2cm  \forall X\le A\ \  |X^{(\p)}/X_{(\p)}|<\infty. \hskip 5cm $$
If $\p$ is either inertial or LIN, to show $(fs)$, we may assume $X_{(\p)}=0$. Thus, since $\p$ is multiplication on  $D_{1}:=Div(A)$, (see Lemma \ref{LemmaDiv}), we have $D_{1}\cap X=0$ and $X$
is reduced. Moreover, by the elementary abelian case, $\p$ is
multiplication on a subgroup of finite index of $A[p]$ and we get that $X[p]$ is finite. It
follows that $X$ is finite. Then (fs) holds by Lemma \ref{Robinson++}.

\medskip\noindent$\bullet$ \emph{ If $\p$ is inertial but not FM, then $A$ is critical,
that is $D:=Div(A)$ has Min and $A/D$ is infinite but bounded.}
 This will be proved by the following steps.

\medskip\noindent I. \emph{\ $A$ is not residually finite}. Assume by contradiction it is. As above, we note that if $\p$ is inertial, then \emph{there is no sequence of subgroups} $X_i$ \emph{with the property} that if we denote $Y_i:=X_i\cap \p(X_i)$ then we have:\\
$(1)$\ \  $Y_{i+1}\cap X_i=Y_i$ \\
$(2)$\ \ the sequence $|\p(X_i)/Y_i|$ is\emph{ strictly increasing}.\\
Otherwise there would exists a
subgroup $X_\omega:=\cup_i X_i$ with the properties that $|\p(X_\omega)/X_\omega\cap \p(X_\omega)|\ge
|\p(X_i)/Y_i|\ge i$ for each $i$.

 On the other hand, we will
construct  a prohibited sequence $X_i$, getting the wished contradiction.
Let $X$ be any finite subgroup of $A$.
By Lemma \ref{Robinson++} the subgroup $K:=X^{(\p)}$ is finite. Since $A$ is residually finite, by  (fs)
there is a $\p$-subgroup $A_{*}$ with finite index in $A$
such that $A_{*}\cap K=0$. Now, as $\p$ is not multiplication on $(A_{*}+K)/K$, there is $a\in A_{*}$ such that $\p(a)\not\in \la a,K\ra$.
Let $Y:=X\cap\p(X)$, $X':=\la a\ra +X$ and $Y':=X'\cap \p(X')$. Let us check that\\
$(1)$\ \  $\p(X)\cap Y'=Y$;\\
$(2')$\ \  $ \p(X')>\p(X)+Y'$.\\
 In fact, to prove $(1)$, if $\p(x)\in \p(X)\cap Y'$ (where $x\in X$),  then $\p(x)=ma+x_0$ with  $m\in\Z,\ x_0\in X$ and $ma=\p(x)-x_0\in A_{*}\cap K=0$, hence $\p(x)=x_0\in Y$ and $(1)$ holds.

To prove $(2')$,
note that if
 $y'\in Y'=X'\cap \p(X')$, then $\exists n,m\in\Z$, $\exists x,x_0\in X$ such
 that $y'=ma+x=n\p(a)+\p(x_0)$. Then $ma-n\p(a)\in A^*\cap K=0$. Hence $x=\p(x_0)\in Y:=X\cap \p(X)$. Since $\p(a)\not\in \la a,K\ra$, then $p$ divides $s$, and so  $Y'\le \la p\p(a)\ra+ Y\not\ni \p(a)$. Therefore $(2')$ holds as $\p(a)\in \p(X')\setminus \la p\p(a)\ra+ Y$.

 Thus we can define by induction a prohibited sequence as above, since from $(1)$ and $(2')$ it follows $|\p(X')/Y'|>|\p(X)/Y|$ and we get  a contradiction.

\medskip\noindent II.
 \emph{\ $A$ is not reduced}.
Otherwise, let $R$ be a basic subgroup of $A$. By  (fs),  $R^{(\p)}/R$ is finite and so $H:=R^{(\p)}$ is residually finite as well. Also,  $A/H$ is divisible. By the above there are  a $p$-adic $\alpha$ and  a finite $\p$-invariant subgroup $A_1$ of $H$ such that $\p=\alpha$ on $H/A_1$. As the kernel $K/A_1$ of $(\p-\alpha)_{|A/A_1}$ contains $H/A_1$ and its image is reduced, while $A/H$ is divisible, it is clear that $K=A$ and $\p={\alpha}$ on $A/A_1$, the  wished contradiction.

\medskip\noindent III.
 \emph{\ $A$ is critical}. Let $A=B\oplus D$ with $B$ infinite but reduced,
$D=Div(A)$ divisible. As (fs) holds at the expense of substituting $A$ with
a finite index $\p$-invariant subgroup $A_0$ we may assume that $B$ and $D$ are both $\p$-invariant. Further, by the above and Lemma \ref{LemmaDiv}, $\p$ is multiplication on both $D$
 and on a finite index subgroup of $B$. So we
may also assume $\p$ is multiplication on $B$. Let $\p$ act on $B$ and $D$ by means of
$p$-adics $\alpha_{1},\alpha_{2}$, resp. As we assumed $\p$ is not FM, we have $\alpha_{1}$ and $\alpha_{2}$ act differently on $B$.

 \emph{If by contradiction $B$ is unbounded}, then there is a
quotient $B/S\iso\Z(p^\infty)$.
 By (fs) we can assume $S$ to be $\p$-invariant and consider $\bar A:=A/S$. This a divisible group
 on which $\p$ acts as a (universal) multiplication by Lemma \ref{LemmaDiv}, contradicting the assumption on $\alpha_{1}$ and $\alpha_{2}$. So that \emph{$B$ is bounded}.

\emph{If by contradiction $D$ has infinite rank}, we may substitute $B$ by $B[p^e]$ where $e$ is the smallest such that $B/B[p^e]$ is finite.
By the reduced case above, $\p$ is multiplication on a subgroup $A_*$ of finite index of $A[p^{\e}]$. Then if $D$ has infinite rank, $\alpha_{1}\equiv\alpha_{2}$ mod $p^\e$ and $\p$ is multiplication on $(B\cap A_*)\oplus D$ which has finite index in $A$, a contradiction.

\medskip\noindent$\bullet$
\ \emph{For finitely many inertial endomorphisms $\p_i$ there are subgroups $A_0$, $B$, $C$, $D$ as in the statement}. Note that if all $\p_i$'s are FM, then clearly there is finite index subgroup $C$ of $A$ such that all $\p_i$'s are multiplication on $C$. Otherwise $A=B_0\oplus D$ is critical, with $B_0$ bounded and $D$ divisible with finite rank. As above, for each $i$ there is a finite index subgroup $B_i$ of $B_0$ such that $\p_i$ is multiplication on $B_i$. Let $B:=\cap_iB_i$. Then $A_0:=B+D$ is the wished subgroup as in the statement.

\medskip\noindent$\bullet$
\emph{\ (FS) holds, if (i), (ii) and (iii) hold}. If each $\p_i$ is FM, then it is multiplication (by the $p$-adic $\alpha_i$) on a subgroup with finite index $A_i$ and modulo a subgroup with finite order $F_i:=im(\p_i-\alpha_i)$. Then all $\p_i$ are multiplication on the finite index subgroup $A_0=\cap_i A_i$ and modulo the finite subgroup $F_0:=F_1\oplus\ldots\oplus F_l$. Therefore for each $X\le A$ we have that $X\cap  A_0$ and $X+F$ are $\p_i$-invariant for each $i$.

If some $\p_i$ is not FM, then $A$ is critical. By the above, there is a finite index subgroup $A_0=B\oplus D$ of $A$ such that each $\p_i$ is multiplication on $B$ and $D$, where $B$ is bounded and $D:=Div(A)$ divisible with finite rank.

Let $X_0:=X\cap A_0$. Then $|X/X_0|\le |A/A_0|$ is finite.
Further, $X_*:=(D\cap X)+(B\cap X)$ is $\p_i$-invariant and the group $X_0/X_*$ is bounded as $B$ is. Also $X_0/X_*$ has rank $r$ at most the
rank of $D$, hence $X_0/X_*$ is finite. Thus
$X/X_{\Phi}$ is finite. Hence each $\p_i$ is inertial.

To show that $X^{\Phi}/X$ is finite as well, we can assume $X_{\Phi}=0$ and that $X$ is finite.
Hence there is $m$ such that $X\le A[m]$. By the above each $\p_i$ is FM on $A[m]$, which is not critical, and $X^{\Phi}$ is finite.

\medskip\noindent$\bullet$
\emph{\ (L) holds}. If $A$ is a $p$-group, notice that we have already seen that LIN implies (fs) and therefore RIN. Thus there are subgroups $A_0,\ B,\ D,\ C$ as in part (R) of the statement.
If $B\ne 0$ (hence $C=0$) we can assume $B$ is infinite and, by Proposition \ref{mult+iner+}.L1, each  $\p_i$ is invertible on $B$ and $\p_i\ne0$ on $D$ (if $D\ne 0$). On the other hand, if some $\p_i$ is not invertible on $C\ne0$ (hence $B\oplus D=0$), then $C$ has Min and we can put $A_0:=Div(C)$ and the statement holds.

In the general case apply Proposition \ref{riduzione_ai_p-gruppi} and deduce that LIN implies RIN since this is true on the $p$-components. Suppose each $\p_i$ is LIN. Then the set  $\pi$ of primes such that some $\p_i$ is not an invertible multiplication on $A_{p}$ is finite, by Proposition  \ref{riduzione_ai_p-gruppi} again.
Then for any $p\in\pi$, if the $p$-component $A_{p}$ is either non-critical or non-Min, then there is  a subgroup $C_p$ with finite index in $A_p$ on which $\p_i$ is invertible multiplication (see Proposition \ref{mult+iner+}.L1). Otherwise
there is a subgroup with finite index in $A_p$ of the shape $B_{p}\oplus D_{p}$, where $B_{p}$ is bounded, $D_{p}$ is divisible of finite rank, such that each $\p_i$ is invertible multiplication on $B_{p}$ and non-zero multiplication on $D_{p}$ (if $D_{p}\not =0$). Then the statement holds with $B:=\oplus_{p\in\pi}B_{p}$, $D:=\oplus_{p\in\pi}D_{p}$ and $C:=A_{\pi'}\oplus \bigoplus\limits_{p\in\pi}C_{p}$.

Conversely, arguing componentwise, it is enough to show that each $\p_i$ is LIN on $B\oplus D$ as in the statement.
Then for each $X\le B\oplus D$ we have $\p(X\cap B)= X\cap B $ and $\p(X\cap D)$ has finite index in $X\cap D$. Then $X/X_*$ and $X_*/X_*\cap\p(X_*)$ are finite, where $X_*:=(X\cap B)+(X\cap D)$.
\qed

\vskip0.8cm
\section{Inertial endomorphisms of a \emph{non-periodic} abelian group}

In order to prove  Theorem A, we give a detailed description of  inertial and LIN endomorphisms of a non-periodic abelian group. We generalize results from \cite{DR} to the more general settings of endomorphisms.
\vskip0.4cm

\begin{Prop}\label{PropMixENDO}
An endomorphism $\p$ of an abelian non-periodic group $A$ is inertial (resp. LIN) if and only if one of the following holds:

\bigskip\noindent
$(a)$ there is a $\p$-invariant finite index subgroup $A_0$ of $A$ on which
$\p$ acts as multiplication by $m\in\Z$
(resp. as multiplication by $\frac{1}{n}$, with $n\ne0$);

\bigskip\noindent
$(b)$ there are finitely many elements $a_i$ such that:

i) the $\p$-submodule $V=\Z[\p]\la a_1,\dots,a_r\ra$
    is torsion-free as an abelian group and $\p$ induces on $V$ multiplication by $\frac{m}{n}$ where $m$, $n$ are coprime
    integers (resp. $m\ne 0$),

ii) the factor group $A/V$ is torsion and $\p$ induces an inertial (resp. LIN)
    endomorphism on $A/V$,

iii) the $\pi(n)$-component of $A$ is bounded.
\end{Prop}

\noindent Notice that $\p=\frac mn$ on both $A/T$ and $V$. Moreover,
in case (b) we have $V\iso\Q^{(n)}\oplus\ldots\oplus\Q^{(n)}$ (finitely many times).

For the proof of the Proposition we need some preliminary result. In the first Lemma we state some well-known facts.

\begin{Lemma}\label{divisibilita1/n} Let $A_1$ be a subgroup of an abelian group $A$ and $\pi$ a set of primes. \\
1)\ If $A/A_1$ is a $\pi'$-group, then $A$ is $\pi$-divisible iff $A_1$ is $\pi$-divisible.\\
2)\ If $A$ if torsion-free, $A/A_1$ periodic and $A_1$ is $\pi$-divisible then
$A/A_1$ is a $\pi'$-group and $A$ is $\pi$-divisible.\\
3)\ If $A_1$ is torsion-free and $\pi$-divisible  while $A/A_1$ is $\pi'$-group, then multiplication by a $\pi$-number is invertible.\qed
\end{Lemma}

\noindent
Next Lemma is a generalization of Lemma \ref{Robinson++} to non-periodic groups.

\begin{Lemma}\label{Robinson+++} Let $A$ be an abelian group, $a\in A$  and $\p\in E(A)$.
If $\p$ is either RIN or LIN, then the torsion subgroup $T$
of $\Z[\p](a)=:\la a\ra^{(\p)}$ is finite. \\
\end{Lemma}
\vskip-0.3cm
\pf We may assume $A=\la a\ra^{(\p)}$. If $a$ has finite order, apply Lemma \ref{Robinson++}.(1). Assume $a$ is aperiodic.
By Proposition \ref{CasoTorsionFree++}, $\p=\frac{m}{n}$ on $A/T$ ($m,n$ coprime), that is
$(n\p-m)(a)$ is periodic. Regard $A$ as $\Z[x]$-module (where $x$ acts as $\p$) and consider
the natural epimorphism  mapping $1$ to $a$ and $x$ to $\p(a)$:
$$F:\Z[x]\mapsto A.$$
 Let $I$ be the inverse image of $T$ via $F$. Then $(nx-m)\subseteq I$ and $\Z[x]/I\iso A/T$ is torsion-free (as $\Z$-module).
Since proper quotients of $\Z[x]/(nx-m)\iso \Z[1/n]=\Q^{(n)}$ are periodic, then $I=(nx-m)$. Applying $F$ we get that
$T=\Z[\p]\la(m\p-n)(a)\ra$ is a cyclic $\p$-submodule. It is finite by Lemma \ref{Robinson++}.
\qed

\medskip
\noindent \textbf{Proof of Proposition \ref{PropMixENDO}}. Assume that $\p$ is either RIN or LIN. By Proposition \ref{CasoTorsionFree++},
 $\p$ is multiplication by $\frac{m}{n}\in\Q$ on $A/T$, where $T:=T(A)$ and $m,n$ are coprime. Let $\pi:=\pi(n)$. We  proceed by a sequence of claims.

 \medskip\noindent$\bullet$
\ \emph{There is a free abelian $F\le A$ such that $V:=\Z[\p]F$ is torsion-free and
$A/V$ is periodic.} In fact,
by Zorn's Lemma, there is a subset $S$ of $A$ which is maximal with respect to
``$F:=\la S \ra$ is free abelian on $S$ and $V:=\Z[\p]\la S \ra$ is torsion-free''. It follows
that $A/V$ is periodic. If not, there is an aperiodic $a\in A$ such that $\la a\ra \cap V=0$.
By Lemma \ref{Robinson+++}, the torsion subgroup of $\Z[\p]\la a \ra$ has finite order $s$. Thus $\Z[\p]\la a^s\ra\iso \Q^{(n)}$
 has rank $1$ (see Proposition \ref{CasoTorsionFree++}) and
$\{ a^s\} \cup S$ has the above properties instead of $S$, a contradiction. Then $\p=\frac mn$ on the torsion-free subgroup $V$ and $A/V$ is periodic.

\medskip\noindent$\bullet$
\ \emph{the $\pi$-component $A_\pi$ is bounded}.
To establish this, assume by contradiction that
 $T$ has a quotient $T/K$ isomorphic to a Pr\"ufer $p$-group,
with $p\in \pi$. By (FS) of Proposition \ref{TEOREMA_PERIODICO+}, $K^{(\p)}/K$ is finite.
Thus, without loss of generality, we may assume $K^{(\p)}=0$,
that is $T$  is a Pr\"ufer
$p$-group. Since $V$ contains a $\p$-invariant subgroup isomorphic to $\Q^{\pi}$, then $T+V$ contains a $\p$-invariant subgroup isomorphic to $\Z(p^{\infty})\oplus\Q^{\pi}$ and it is enough to check that:\\  \emph{- the endomorphism
 $\p=\alpha\oplus \frac{m}{n}$ of $\Z(p^{\infty})\oplus\Q^{\pi}$  is neither RIN nor LIN} ($\alpha$ a $p$-adic). To this aim consider
the ``diagonal" subgroup $$H=\{\ [tp^{-i}]\oplus tp^{-i}\ |\ i\in\N, t\in\Z, \ [tp^{-i}]\in \Q^{(p)}/\Z \}$$ and note $(m-n\p)([p^{-i}]\oplus p^{-i})=(m-n\alpha)[p^{-i}]\oplus 0.$
Since $p$ does not divide $m$, we have $\Z(p^\infty)= (m-n\p)(H)$.
Therefore $H+\p(H)$ is commensurable neither to $H$ nor to $\p(H)$, the wished contradiction.

\medskip\noindent$\bullet$
\ \emph{If $r_{0}(A)$ is finite, then $(b)$ holds}. This is now clear.

\medskip\noindent$\bullet$
\ \emph{If $A$ has not FTFR and $\p$ is inertial, then $(a)$ holds}.
By Proposition \ref{CasoTorsionFree++} we
know that $\p=m\in\Z$ on $V=F$ as above, which is a free abelian group on the basis $S$. So
there exists a ($\p$-invariant) subgroup $W$
such that $V/W$ is a periodic group whose $p$-component is divisible
with infinite rank for each prime $p$.
By the torsion case, Proposition \ref{TEOREMA_PERIODICO+}, $\p$ is FM
on $A/W$, since it contains a divisible $p$-group with infinite rank for each prime $p$.
Without loss of generality, we can assume $\p$ is a multiplication indeed on
$A/W$ and, for each $p$, $\p$ is represented by the $p$-adic $\alpha_p$ on the $p$-component of $A/W$.
Then $\alpha_p=m$, as they must agree on the $p$-component of $V/W$, which it is unbounded.
Then  $W\ge im(\p-m)\iso {A}/{ker(\p-m)}$ where the former is torsion-free and
the latter is periodic, as a factor of $A/V$. Thus $\p=m$ on $A$.

\medskip\noindent$\bullet$
\ \emph{If $A$ has not FTFR and $\p$ is LIN, then $(a)$ holds}.
Let  $F$ and $V=\Z[\p]F$ as in the first claim above. By Proposition \ref{CasoTorsionFree++}, $m=1$, that is $\p=\frac{1}{n}$ on $V$. Then $V/F$ is the sum of infinitely many copies of $\Z(p^\infty)$ for each prime $p\in\pi:=\pi(n)$. Take $F_*$ such that  $F/F_*$ is a periodic $\pi'$-group whose $p$-component is divisible with infinite rank for each prime $p\in\pi'$. Let $V_*/F_*$ be the $\pi$-component of $V/F_*$. Since
$V/V_*$ is a $\pi'$-group, $V_*$ is $\pi$-divisible, by Lemma \ref{divisibilita1/n}.(1). Thus $V_*\le V$ is $\p$-invariant and $V/V_*\iso_\p F/F_*$.
Let $A_*/V_*$ be the $\pi'$-component of $A/V_*$.

\emph{ We claim $A/A_*$ is finite}. It is enough to check that \emph{the $\pi$-component $A_1/V_*$ of $A/V_*$ is finite}. To this aim, notice that $T_1:=T(A_1)=A_{\pi}$. On one hand $A_1/(T_1+V_*)$ is a $\pi$-group by definition of $A_1/V_*$; on the other hand
$A_1/(T_1+V_*)$ is $\pi'$-group as $(T_1+V_*)/T_1$ is $\pi$-divisible and $A_1/T_1$ is torsion-free (see Lemma \ref{divisibilita1/n}.(2)). Thus
$A_1=T_1\oplus V_*$ and the claim reduces to show $T_1$ is finite.

\emph{Assume by contradiction that $T_1$ has infinite rank} (since by $(iii)$ above, $T_1$ is  bounded).
By Proposition \ref{TEOREMA_PERIODICO+}, $\p$ is FM on $T_1$.
Then $\p$ is a multiplication  by an integer $s$ not multiple of $p$
on a countable $\Z_p$-submodule $B=\oplus_i \la b_i\ra\le T_1$. Let $\{a_{i}|\,i<\omega\}$ be a countable subset of the above basis $S$ for $F$ and set $W:=\la V_i\ |\ i<\omega\ra = \oplus_i V_i$, where $V_i:=\Z[\p]a_i$. Also,
let $M:=B \oplus W$ and $H:=\la a_i+b_i\ |\ i<\omega\ra$ its ``diagonal" subgroup, which is free on the $\Z$-basis of the $a_i+b_i$'s.
Since $\p$ is one-to-one on $M$, then $\p(H)$ is torsion-free as $H$ is.
For all $i$ we have:
$H+\p(H)\ni (p-n\p)(a_i+b_i)=(p-1)a_i,$ as $p$ divides $n$. Since  $pa_{i}\in H$, then $a_i\in H+\p(H)$. Thus
$B\le H+\p(H)$. Therefore $(H+\p(H))/\p(H)\ge (B+\p(H))/\p(H)\iso B$ is infinite, contradicting $\p$ is LIN. \emph{Thus  $A/A_*$ is finite}.

Let us show that \emph{$\p=\frac{1}{n}$ on some $A_0$ with finite index in $A_*$.}
As  $A_*/V_*$ is a $\pi'$-group and, for each prime $p\in\pi'$,
its $p$-component contains a divisible $p$-group of infinite rank, by Proposition \ref{TEOREMA_PERIODICO+} we have that
$\p$ is FM on $A_*/V_*$. Thus $\p$ is multiplication on some $A_0/V_*$ with finite index in $A_*/V_*$.
On one side $\p=\frac{1}{n}$ on $V$. On the other side,
by Lemma \ref{divisibilita1/n}.(3) the multiplication by $\frac{1}{n}$ is an endomorphism on the whole $A_{0}$. Then, as $ker(\p_{|A_0}-\frac{1}{n})\le V$ and $(\p-\frac{1}{n})(A_0)\le V_*$, we have that
 $A_0/ ker(\p_{|A_0}-\frac{1}{n})\iso (\p-\frac{1}{n})(A_0)$ is both periodic and torsion-free. Therefore $\p=\frac{1}{n}$ on $A_0$. Thus $(a)$ holds if $r_{0}(A)$ is infinite.


\bigskip

\emph{Conversely}, if $\p$ is as in $(a)$, it is trivial that $\p$ is RIN (or LIN, resp.).
 Let then $\p$ as in $(b)$. We have to show that for each subgroup $X$ of $A$ the statement
$R(X)$ (resp. $L(X)$) below holds.
$$ R(X):=\ \left(\ \left|\dfrac{X+\p(X)}{X}\right|<\infty\ \right)\ \ \ \ \ \
\ L(X):=\ \left(\ \left|\dfrac{X+\p(X)}{\p(X)}\right|<\infty\ \right)$$
Let $\pi:=\pi(n)$. We proceed by steps:
\medskip

\medskip\noindent$\bullet$  \emph{$\p=\frac{m}{n}$ is inertial on $A/T$ and so $A/T$ is $\pi$-divisible}.
 In fact, if $a\in A$, there is a non-zero integer $s$ such that $sa\in V$.
 Thus $s(n\p-m)(a)=(n\p-m)(sa)=0$ and $(n\p-m)(A)\subseteq T$, as claimed.

\medskip

\medskip\noindent$\bullet$  \emph{If $X$ is periodic, then $R(X)$ (resp. $L(X)$) holds}. This is very easy, since $X^{(\p)}\cap V=0$ and one can verify $R(X)$ (resp. $L(X)$) mod $V$.

\medskip
\medskip\noindent$\bullet$ \emph{It is enough to show that $R(X_0)$ (resp $L(X_0)$) holds for each torsion-free subgroup $X_0$ of $A$}.
Let $X$ be any subgroup and $U=T(X)$. By the above $\p$ induces on $T$ a RIN (resp LIN) endomorphism. By Proposition \ref{TEOREMA_PERIODICO+} applied to $T$, we have (FS), so that
$U/U_{(\p)}$ is finite. Since the hypotheses hold modulo $U_{(\p)}$ (which is periodic), that is
for the  endomorphism induced by $\p$ on the group $A/U_{(\p)}$, we can assume $U_{(\p)}=0$ that is $U=T(X)$ is finite. Therefore $X=X_0\oplus U$ splits on $U$. Since $X/X_0$ is finite, then $R(X_0)$ (resp. $L(X_0)$) implies
straightforward  $R(X)$  (resp. $L(X)$).

\medskip

\medskip\noindent$\bullet$ \emph{We may assume $A_\pi=0$}. Recall that by hypotesis $B:=A_\pi$ is bounded by some $e$.
 Let $X\le A$ be torsion-free. Clearly $X$ has finite rank. If $\p$ is RIN on $A/B$, let
$\left|\frac{X+\p(X)+B}{X+B}\right|=s<\infty$. Then $s\p(X)\le X+B$. Thus $es\p(X)\le X$. It follows $\p(X)+X/X$ is finite, as claimed, as it is bounded and has finite rank.  If $\p$ is LIN on $A/B$, let $\left|\frac{X+\p(X)+B}{\p(X)+B}\right|=s<\infty$. By an argument as above, $esX\le \p(X)$ and $(\p(X)+X)/\p(X)$ is finite.

\medskip

\medskip\noindent$\bullet$  \emph{If $A_\pi=0$ and $X$ is torsion-free
then R(X) (resp. L(X)) holds.}
As the hypotheses on $\p$ hold even in
$A_1:=\Z[\p]X$ with respect to $V_1:=A_1\cap V$, we can assume $A=A_1$, that is
$X$ has maximal torsion-free rank $r$ and $V\iso \Q^{\pi}\oplus\ldots\oplus \Q^{\pi}$ ($r$ times).

 Let $K/X$ be the $\pi$-component of $A/X$ (which is periodic). By hypothesis $R(K+V)$ holds, thus $R(K)$ holds, since  $(K+V)/K\iso V/(V\cap K)$ is finite as it is a $\pi'$-group.
On the other hand, $K$ is  torsion-free, as $T$ is a $\pi'$-group. Thus $T(K+\p(K))$ is finite.

Let $Y:=X+\p(X)$, $Y_R:= Y\cap (X+T)$ and $Y_L:=Y\cap (\p(X)+T)$.
On one side, $Y_R/X$ and $Y_L/\p(X)$ are both finite, as  isomorphic to quotients of $Y\cap T\le T(K+\p(K))$, which is finite.
On the other side, $|Y/Y_R|=|(Y+T)/(X+T)|<\infty$ because of $R(X+T)$ (resp. $|Y/Y_L|=|(Y+T)/(\p(X)+T)|<\infty$ because of $L(X+T)$)
 which holds as $\p$ is inertial on $A/T$ (and $m\ne0$, resp.) see Proposition \ref{CasoTorsionFree++}.
\QED

\section{Proofs of the main results}

We use a Lemma dealing with finitely many inertial endomorphisms. Denote by $\oplus_r \Q^{\pi}$ the direct sum of $r$ copies of $\Q^{\pi}$.

\begin{Lemma}\label{TeomixENDOmany}
Let $\p_1,\ldots,\p_t$ be finitely many inertial endomorphisms of an abelian group $A$ with FTFR, where $\p_i=\frac{m_i}{n_i}\in\Q$ on $A/T$ ($m_{i}$, $n_{i}$ coprime $\forall i$) and $\pi:=\pi(n_1\cdot...\cdot n_l)$. Then there is a $\Z[\p_1,...,\p_t]$-submodule $V\iso \oplus_r \Q^{\pi}$ ($r\in\N$) and such that $A/V$ is periodic.
\end{Lemma}

\pf It is easily seen that for each $p_j\in\pi$ there is $\psi_j\in\Z[\p_1,...,\p_t]$ such that $\psi_j=r_j/p_j$ on $A/T$ ($r_j$ and $p_j$ coprime). Then there are coprime $m,n$ such that $\p:=\sum_j \psi_j= m/ n$ on $A/T$ where $\pi( n)=\pi$.

By Proposition \ref{PropMixENDO} for each $i$ there is a torsion-free $\p_i$-invariant subgroup $V_i$ such that $A/V_i$ is periodic.
Pick $\Z$-independent elements $b_{1},\dots,b_{r}$ of $\cap_i V_{i}$, where $r=r_0(A)$.
By Lemma \ref{Robinson+++}, for each $k$ there exists $a_k\in\la b_k\ra$ such that
$\Z[\p]\la a_k\ra$ is torsion-free with rank $1$. Therefore the subgroup
 $V:=\Z[\p]\la a_{1},\dots,a_{r}\ra$
 is  torsion-free and $\pi$-divisible.

 We claim $\p_i(V)\subseteq V$. Set $W_i:=\Z[\p_i]\la a_{1},\dots,a_{r}\ra$, which is torsion-free as it is a subgroup of $V_{i}$ and is contained in $V$, since $\p_{i}=\frac{m_i}{n_i}$ on $W_i$.  Then $V/W_i$ is a $\pi$-group, as $V/\la a_1,...,a_r\ra$ is such.

 Let $a$ be any element of $V$ and $e$ be the bound of $A_{\pi}$ (which is bounded by Proposition \ref{PropMixENDO}). Therefore there is a $\pi$-number $t$ such that $ta\in W_i$ hence
$(n_i\p_i-m_i)(ta)\in T\cap W_i=0$. Thus $(n_i\p_i-m_i)(a)\in A_\pi$, so that $e(n_i\p_i-m_i)(a)=0$. Since $e$ and $n_i$ are $\pi$-numbers and $V$ is $\pi$-divisible, we have
$\p_i(V)=en_i\p_i(V)=em_iV\subseteq V$ as wished.
\QED

\noindent\textbf{Proof of Theorem A}.
Assume all $\p_{i}$ are inertial. If  $A$ has not FTFR, then by  Proposition \ref{PropMixENDO}
it is clear that $(a)$  holds.

Assume now $A$ has FTFR. As in Lemma \ref{TeomixENDOmany}, there is  a $\Z[\p_1,...,\p_t]$-submodule $V$ such that $V\iso \oplus_r \Q^{\pi}$ ($r\in\N$) and  $A/V$ is periodic. Let $\pi_2$ be the set of primes $p$ such that some $\p_i$ is not FM on the $p$-component of $A/V$.
Note that the definition of $\pi_2$ is independent of $V$, as all possible $V$ are commensurable each other.
From Proposition \ref{TEOREMA_PERIODICO+} it follows that $\pi_2$ is finite.
Let $\pi_1:=\pi\cup \pi_2$. On one side, $A_\pi$ is bounded, by Proposition \ref{PropMixENDO}. On the other side, for each $p\in \pi_2$ the $p$-component $A_p$ of $A$ is the sum of a bounded subgroup and a finite rank divisible subgroup. Thus there is $C^*$ such that $A=A_{\pi_1}\oplus C^{*}.$

By  Proposition \ref{TEOREMA_PERIODICO+} and the defintion of $\pi_{1}$, there is a finite index subgroup $B\oplus D$ of $A_{\pi_1}$ such that $B$ is bounded, $D$ is divisible with finite rank hence a $\pi'$-group and each $\p_i$ acts by multiplications on both $B$ and  $D$, as in the statement.

{We may assume $V\le C^{*}$}. In fact $|V/(V\cap C^{*})|=:s$ is finite as $V\iso \oplus_r \Q^{\pi}$ ($r\in\N$) and $V/(V\cap C^{*})\iso (V+C^*)/C^{*}$ is periodic with bounded $\pi$-component. So we may  substitute $V$ with $sV$ and we get $V\le C^{*}$.

Use bar notation in $\bar A:=A/V$. Consider the primary decomposition $\bar C^{*}=\bar C_1^{*}\oplus \bar C_{0}^{*}$ where $\bar C_1^{*}$ (resp. $\bar C_0^{*}$) is a $\pi_1$-group (resp. $\pi_1'$-group). By  (FS) of Proposition \ref{TEOREMA_PERIODICO+} and the definition of $\pi_1$, each $\p_i$ is multiplication on a subgroup $\bar C_{0}$ with finite index in $\bar C_{0}^{*}$. On the other side $C_{1}^{*}$ is torsion-free (with finite rank) hence $\bar C_1^{*}$ has Min, thus $\bar C_{1}^{*}$
has a divisible finite index subgroup, say $\bar C_{1}$, on which each $\p_i$ is multiplication (see Lemma \ref{LemmaDiv}).
Therefore there is a finite index subgroup $C:=C_{1}+ C_{0}$ of $C^*$ such that  $\p_i$ is multiplication on $\bar C$. So conditions $(i)$, $(ii)$, $(iii)$ of the statement hold for $A_{0}:=B\oplus D\oplus C$.

\emph{To prove condition (iv)}, for each $i$ let $\p_i=\frac{m_{i}}{n_{i}}\in\Q$ on $V$ (see Proposition \ref{CasoTorsionFree++}) and $p\in\pi(D)$ such that the $p$-component $\bar C_p$ of $\bar C$ is infinite (hence unbounded).
On one hand, as $C_{p}$ is torsion-free, we get that $\p_i=\frac{m_{i}}{n_{i}}$ on $C_p$. On the other hand, $\p_i$ acts by the same $p$-adic $\alpha$ on $D_p$ as on $\bar D_p$. Therefore $\alpha=\frac{m_{i}}{n_{i}}$ by Proposition \ref{TEOREMA_PERIODICO+} (as $\bar D_p\oplus \bar C_p$ is non-critical).

\medskip
Conversely, for the sufficency of the condition, it is clear that  if $(a)$ holds, then all $\p_i$ are inertial, so  only case $(b)$ is left. By Proposition \ref{inerXfin} we may assume $A=A_0$. Fix $i$ and $a_1,...,a_r$ such that $V:=\Z[\p_1,...,\p_t]\la a_1,....a_r\ra$. Let $V_i:=\Z[\p_i]\la a_1,....a_r\ra$. Then $V/V_i$ is a divisible $\pi$-group with finite rank. By  Proposition \ref{PropMixENDO}, any $\p_i$ is inertial iff it is such on the periodic group $A/V_i$. Clearly, by  Proposition \ref{TEOREMA_PERIODICO+}, $\p_i$ is already inertial on $A/V$. Use bar notation  in $\bar A:=A/V_i$. By Proposition \ref{riduzione_ai_p-gruppi}, it is enough to show that \emph{ $\p_i$ is inertial on all  $p$-components of $\bar A$ and even multiplication on all but finitely many}.

If $p\in\pi$ (which is finite), the $p$-component of $\bar{A}$ is $\bar A_p\oplus \bar {C}_p$ where $A_p$ is the $p$-component of $A$ and $\bar C_p$, the $p$-component of $\bar C$, contains $\bar V_p$, the $p$-component of $\bar  V$. Since $C$ has no elements of order $p$, then $C_p$ is torsion-free and, by Lemma \ref{divisibilita1/n}, $\bar {C}_p/\bar V_p$ is a $\pi'$-group.
Hence $\bar {C}_p=\bar V_p$ is  divisible of finite rank. Therefore $\p_i$ is multiplication on $\bar {C}_p$. On the other hand,  $\p_i$ is multiplication even on $\bar A_p\iso_{\p_i} A_p\le B$, which is bounded. Thus $\p_i$ is inertial on $\bar A_p\oplus \bar {C}_p$ by Proposition \ref{TEOREMA_PERIODICO+}.

If $p\not\in\pi$, then the $p$-component of $\bar A$ is $\p_i$-isomorphic to the $p$-component of $A/V$ as $V/V_{i}$ is a $\pi$-group.
\qed

\medskip
\noindent \textbf{Proof of  Corollary A.} Apply Theorem A to any couple of inertial endomorphisms $\p_1$, $\p_2$ of an abelian group $A$. If $(a)$ holds for both $\p_1$ and $\p_2$, then $\p_1\p_2-\p_2\p_1=0$ on a subgroup with finite index of $A$, since multiplications commute. Otherwise, by Proposition \ref{CasoTorsionFree++}, $\p_1$, $\p_2$
 commute on $A/T$ anyway, where $T=T(A)$. Moreover, there is a subgroup $A_0$ with finite index of $A$ such that $\p_1$, $\p_2$ commute on $A_0/V$ for some $V$ as in Theorem A. As $T\cap V=0$, then $\p_1$, $\p_2$ commute on $A_0$, as wished.\qed

\medskip
\noindent \textbf{Proof of  Theorem B.} It is enough to prove the statement for a finitely generated subgroup $\G$ of $IAut(A)$. Let $A_0$ and $V$ as in Theorem A and $T_0:=T(A_0)$. Then $\G'$ acts trivially on both $A_0/V$ and $A_0/T_0$. Thus $\G'$ acts trivially on $A_0$ and (1) holds.

The subgroup $\G_0:=C_\G(A/A_0)$ has finite index in $\G$. On the other hand, by the above, $\Sigma:=\G'\cap \G_0$ stabilizes the series $0\le A_0\le A$ and embeds in $Hom(A/A_0,A_0)$, which is bounded by $m:=|A/A_0|$. Thus $[A,\Sigma]\le A_0[m]=B[m]\oplus C[m]\oplus D[m]$. As each $\g\in\G$ acts by multiplications on $B[m], C[m], D[m]$ and $\G$ is finitely generated, $|\G/\G_2|$ is finite, where $\G_2=C_{\G_0}([A,\Sigma])$. On the other hand, we have that $[A,\Sigma,\G_2]=0$ and $[A,\G_2,\Sigma]\le[A_0,\Sigma]=0$. Thus, by the Three Subgroup Lemma $[\G_2,\Sigma,A]=0$, that is
$\Sigma$ is contained in the center of $\G_2$ which turns to be nilpotent and finitely generated as well. Therefore $\G$ has the maximal condition (Max) and $\G'$ is finite, being periodic. Finally, as $\G$ is finitely generated, we have $\G$ is central-by-finite.\qed

\medskip
\noindent \textbf{Proof of  Corollary B.} In part $(1),$ by Proposition \ref{TEOREMA_PERIODICO+} and \ref{PropMixENDO}, LIN implies RIN. Also, in any case, if $\p$ is invertible,  then $\p$ has LIN iff the $\p^{-1}$ has RIN.

For statement $(2)$, note that if $\g_1$ is RIN and $\g_2$ is LIN, then $\g_1\g_2=\g_2\g_1[\g_1,\g_2]$ where
$\g_1[\g_1,\g_2]$ is RIN as $[\g_1,\g_2]$ is finitary.
\qed

\medskip
\noindent \textbf{Proof of  Proposition A.} Let  $A$ be any abelian group.
If $V$ is a torsion-free and $A/V$ is periodic, denote $T:=T(A)$ and $\bar A:=A/(V+T)$, then the stabilizer $\Sigma$ of the series $0\le (V+T)\le A$ is canonically isomorphic to $Hom(\bar A,V+T)=Hom(\bar A,T)$.  In the particular case that $V$ has finite rank and $A/V$ is locally cyclic, then by Proposition \ref{PropMixENDO} we have
$\Sigma\le IAut(A)$ and $\Sigma\cap FAut(A)$ corresponds to the subgroup of $Hom(\bar A,T)$ formed by the homomorphisms with finite image. Therefore, if the abelian group $A$ is such that, for each prime $p$, $A_p$ has order $p$ while the $p$-component of $A/V$ has order $p^2$, we have:  $\Sigma\iso Hom(\bar A,T)\iso\prod_p \Z(p)$ and $ \Sigma\cap FAut(A)\iso \bigoplus_p \Z(p)$, where  $p$ ranges over the set of all primes.

To show the existence of the wished group $A$, let $G:=B\oplus C$ where
$B:=\prod_p\langle b_{p}\rangle$,
$C:=\prod_p\langle c_{p}\rangle$, and $b_p$, $c_p$ have order $p$, $p^2$ (resp.) .  Consider the (aperiodic) element $v:=(b_p+pc_p)_p\in G$ and $V:=\la v\ra$. Note that for each prime $p$ there exists an element $d_{(p)}\in G$ such that $pd_{(p)}=v-b_{p}$. Define
$A:=V+\langle d_{(p)}|\ p\ \rangle$.
 Then we have that $A/T\iso \la 1/p\ |\ p\ \ra\le \Q$ as it has torsion free rank $1$ and $v+T$ has $p$-height $1$ for each prime. Then $T=T(B)\iso\bigoplus_p \Z(p)$, while
 the $p$-component of $A/V$ is generated by $d_{(p)}+V$ and has order $p^2$ as $pd_{(p)}=v-b_{p}$.
\qed


{
}


\end {document}